\documentclass[10pt]{amsart} 
\usepackage{pstricks, pst-node, amssymb}
\usepackage{epsfig,color}
\usepackage{graphicx}
\usepackage{amsmath, amsfonts}
\usepackage{amsthm,amssymb,amsbsy, latexsym, stmaryrd, amscd, xy} 
\usepackage[mathscr]{eucal}
\usepackage{t1enc}
\usepackage{fancyheadings} 
\usepackage{epic} 
\usepackage{eepic}

\psset{unit=1.1pt, arrowsize=4pt, linewidth=.7pt}
\psset{linecolor=blue}
\newgray{grayish}{.90}
\newrgbcolor{embgreen}{0 .5 0}
\def\vblack(#1, #2)#3{\cnode*[linecolor=black](#1, #2){3}{#3}}
\def\vwhite(#1,#2)#3{\cnode[linecolor=black,fillcolor=white,fillstyle=solid](#1,#2){3}{#3}}
\countdef\x=23
\countdef\y=24
\countdef\z=25
\countdef\t=26

\def\tbox(#1,#2)#3{
\x=#1 \y=#2
\multiply\x by 12
\multiply\y by 12
\z=\x \t=\y
\advance\z by 12
\advance\t by 12
\psline(\x,\y)(\x,\t)(\z,\t)(\z,\y)(\x,\y)
\advance\x by 6
\advance\y by 6
\rput(\x,\y){{\bf #3}}}

\newtheorem{prop}{Proposition} 
\newtheorem{lemma}[prop]{Lemma}

\newtheorem{corollary}[prop]{Corollary} 
\newtheorem{theorem}[prop]{Theorem}
\newtheorem{thm}[prop]{Theorem}

\theoremstyle{definition} 
\newtheorem{definition}[prop]{Definition}
 
\newtheorem{example}[prop]{Example}

\newcommand{\pref}[1]{{(\protect\ref{#1})}}
\def\dd{\kern.4ex\mbox{\raise.7ex\hbox{{\rule{.45em}{.12ex}}}}\kern.4ex}

\def\bxcxa{\ensuremath{(2{\dd}3\dd1)}}
\def\axcb{\ensuremath{(1{\dd}32)}}
\def\acxb{\ensuremath{(13{\dd}2)}}
\def\cbxa{\ensuremath{(32{\dd}1)}}
\def\bxac{\ensuremath{(2{\dd}13)}}

\def\abc{\ensuremath{(123)}}

\def\baxc{\ensuremath{(21{\dd}3)}}
\def\cxba{\ensuremath{(3{\dd}21)}}

\def\bxca{\ensuremath{(2{\dd}31)}}
\def\caxb{\ensuremath{(31{\dd}2)}}

\def\newMAH#1{%
\expandafter\def\csname #1\endcsname{\mathop{\mbox{{\sc#1}}}\nolimits}%
}

\newMAH{maj}

\def\newexpMAH#1{%
\expandafter\def\csname exp#1\endcsname{\mathop{\mbox{{\footnotesize\sc{#1}}}}\nolimits}%
}

\newexpMAH{maj}

\def\Le{\hbox{\rotatedown{$\Gamma$}}}

\def\biw#1,#2,{%
\begin{pmatrix}
#1\cr
#2\cr
\end{pmatrix}
}

\def\ch#1,#2,{{#1\choose #2}}

\newcommand{\cls}{\mathcal{S}}

\newcommand{\ra}{\longrightarrow} 

\def\newop#1{\expandafter\def\csname#1\endcsname{\mathop{\rm#1}\nolimits}}

\def\newscop#1{%
\expandafter\def\csname #1\endcsname{\mathop{\mbox{{\sc#1}}}\nolimits}%
}

\def\newcapop#1{%
\expandafter\def\csname #1\endcsname{\mathop{\mbox{{\sc#1}}}\nolimits}%
}

\newscop{inv}
\newscop{mak}
\newscop{makl}

\newcommand\st{\; | \;}
\def\emm#1,{{\em #1}}

\def\bskip#1,{\vspace*{#1\baselineskip}}

\newop{des}
\newop{wex}

\newscop{wb}
\newscop{wt}

\newscop{db}
\newscop{dt}

\newscop{wextopsum}
\newscop{wexbotsum}
\newscop{destopsum}
\newscop{desbotsum}

\newscop{wextopset}
\newscop{wexbotset}
\newscop{destopset}
\newscop{desbotset}

\newscop{rembr}

\def\aa{\mathrm{{A}_{EE}}}
\def\ab{\mathrm{{A}_{NN}}}
\def\ac{\mathrm{{A}_{EN}}}
\def\ad{\mathrm{{A}_{NE}}}

\def\ca{\mathrm{{C}_{EE}}}
\def\cb{\mathrm{{C}_{NN}}}

\newcapop{cab}

\newcommand\Aa{\aa}
\newcommand\Ab{\ab}
\newcommand\Ac{\ac}
\newcommand\Ad{\ad}

\newcommand\Ca{\ca}
\newcommand\Cb{\cb}

\newcommand\T{\ensuremath{\mathcal{T}_n^k}}

\catcode`\@=11
\def\section{\@startsection{section}{1}%
 \z@{.7\linespacing\@plus\linespacing}{.5\linespacing}%
 {\normalfont\bfseries\scshape\centering}}
\def\subsection{\@startsection{subsection}{2}%
  \z@{.5\linespacing\@plus\linespacing}{.5\linespacing}%
  {\normalfont\bfseries\scshape}}

\def\subsubsection{\@startsection{subsubsection}{3}%
  \z@{.5\linespacing\@plus.7\linespacing}{-.5em}%
  {\normalfont\itshape}}
\catcode`\@=12
\title{Permutation tableaux and permutation patterns} 

\author{Einar Steingr\'{\i}msson \and Lauren K. Williams} 

\address{Reykjavik University, Reykjavik, Iceland, and Dept. of
  Mathematics, Chalmers Univ. of Technology, Göteborg, Sweden}
\email{einar@math.chalmers.se} 
\address{
  Department of Mathematics, Harvard, Cambridge, MA 02138}

\email{lauren@math.harvard.edu} 

\thanks{The first author was partially supported by the European
Commission's IHRP Programme, grant HPRN-CT-2001-00272, ``Algebraic
Combinatorics in Europe''}

\date{\today}

\begin{document}

\begin{abstract}
In this paper we introduce and study a class of tableaux which we call
permutation tableaux; these tableaux are naturally in bijection with
permutations, and they are a distinguished subset of the
$\Le$-diagrams of Alex Postnikov \cite{Postnikov, Williams}.  The
structure of these tableaux is in some ways more transparent than the
structure of permutations; therefore we believe that permutation
tableaux will be useful in furthering the understanding of
permutations.  We give two bijections from permutation tableaux to
permutations.  The first bijection carries tableaux statistics to
permutation statistics based on relative sizes of pairs of letters
in a permutation and their places.  We call these statistics \emm weak
excedance statistics, because of their close relation to weak
excedances.  The second bijection carries tableaux statistics (via the
weak excedance statistics) to statistics based on generalized
permutation patterns.  We then give 
enumerative applications of
these bijections.  
One nice consequence of these results is that the polynomial
enumerating permutation tableaux according to their content
generalizes both Carlitz' $q$-analog of the Eulerian numbers 
\cite{Carlitz} and the more recent $q$-analog of the Eulerian numbers
found in \cite{Williams}.
We conclude our paper with a list of open problems, as well as
remarks on progress on these problems which has been made by 
A. Burstein, S. Corteel, N. Eriksen, A. Reifegerste, and X. Viennot.
\end{abstract}

\maketitle \thispagestyle{empty}

\section{Introduction}

The aim of this article is to advertise a new class of tableaux
together with two curious bijections for the study of permutations.
We call these tableaux {\em permutation tableaux}; they are naturally
in bijection with permutations, and are a distinguished subset of Alex
Postnikov's $\Le$-diagrams \cite{Postnikov}, which were enumerated by
the second author \cite{Williams} because of their connection with the
totally nonnegative part of the Grassmannian.

Recall that a {\em partition} $\lambda = (\lambda_1, \dots,
\lambda_k)$ is a weakly decreasing sequence of nonnegative
integers. For a partition $\lambda$, where $\sum \lambda_i = m$, the
{\em Young diagram} $Y_\lambda$ of shape $\lambda$ is a left-justified
diagram of $m$ boxes, with $\lambda_i$ boxes in the $i$-th row.

We define a {\em permutation tableau} $\T$ to be a partition $\lambda$
such that $Y_\lambda$ is contained in a $k \times (n-k)$ rectangle,
together with a filling of the boxes of $Y_\lambda$ with $0$'s and
$1$'s such that the following properties hold:
\begin{enumerate}
\item Each column of the rectangle contains at least one $1$.
\item There is no $0$ which has a $1$ above it in the same column
{\em and} a $1$ to its left in the same row.
\end{enumerate}

We call such a filling a {\em valid} filling of $Y_\lambda$.  Observe
that the requirement in (1) implies that the Young diagram must have $n-k$
columns, whereas the number of rows may be smaller than $k$.

If we forget the requirement  (1) above we recover the definition of
a $\Le$-diagram~\cite{Postnikov}.

Figure \ref{PermTab} gives an example of a permutation tableau.

\begin{figure}[h]
\pspicture(-95,-12)(230,94)
\rput(190,36)
{$\begin{array}{l}
k=7,\ n=17\\
\lambda=(10,9,9,8,5,2)
\end{array}$}
\rput(-10,36){$k$}
\rput(60,82){$n-k$}
\psline[linecolor=black,linewidth=0.5pt]{-}(0,-12)(120,-12)(120,72)(0,72)(0,-12)
\tbox(0,0){1}
\tbox(1,0){1}
\tbox(0,1){0}
\tbox(1,1){0}
\tbox(2,1){0}
\tbox(3,1){1}
\tbox(4,1){1}
\tbox(0,2){0}
\tbox(1,2){0}
\tbox(2,2){0}
\tbox(3,2){0}
\tbox(4,2){0}
\tbox(5,2){0}
\tbox(6,2){1}
\tbox(7,2){1}
\tbox(0,3){0}
\tbox(1,3){0}
\tbox(2,3){0}
\tbox(3,3){0}
\tbox(4,3){0}
\tbox(5,3){0}
\tbox(6,3){0}
\tbox(7,3){0}
\tbox(8,3){0}
\tbox(0,4){1}
\tbox(1,4){1}
\tbox(2,4){1}
\tbox(3,4){1}
\tbox(4,4){0}
\tbox(5,4){1}
\tbox(6,4){1}
\tbox(7,4){1}
\tbox(8,4){1}
\tbox(0,5){0}
\tbox(1,5){1}
\tbox(2,5){1}
\tbox(3,5){0}
\tbox(4,5){0}
\tbox(5,5){1}
\tbox(6,5){0}
\tbox(7,5){1}
\tbox(8,5){0}
\tbox(9,5){1}
\endpspicture
\caption{A permutation tableau}
\label{PermTab}
\end{figure}

Note that there is a unique permutation tableau which has $n$ rows and 
no columns!

We will also think of a permutation tableau $\T$ as a $k \times
(n-k)$ array of $0$'s, $1$'s, and $2$'s, by simply taking the previous
description of a permutation tableaux and putting a $2$ in every box
of the rectangle which is not in $Y_\lambda$, as in
Figure~\ref{PermTab2}.  We position the partition shape so that
its top row lies at the top of the rectangle; therefore the 2's
cut out a (rotated) Young diagram in the southeast corner of the
rectangle.

\begin{figure}[h]
\pspicture(-95,-12)(230,94)
\rput(190,36)
{$\begin{array}{l}
k=7,\ n=17\\
\lambda=(10,9,9,8,5,2)
\end{array}$}
\rput(-10,36){$k$}
\rput(60,82){$n-k$}
\psline[linecolor=black,linewidth=0.5pt]{-}(0,0)(120,0)(120,72)(0,72)(0,0)
\tbox(0,-1){2}
\tbox(1,-1){2}
\tbox(2,-1){2}
\tbox(3,-1){2}
\tbox(4,-1){2}
\tbox(5,-1){2}
\tbox(6,-1){2}
\tbox(7,-1){2}
\tbox(8,-1){2}
\tbox(9,-1){2}

\tbox(0,0){1}
\tbox(1,0){1}
\tbox(2,0){2}
\tbox(3,0){2}
\tbox(4,0){2}
\tbox(5,0){2}
\tbox(6,0){2}
\tbox(7,0){2}
\tbox(8,0){2}
\tbox(9,0){2}

\tbox(0,1){0}
\tbox(1,1){0}
\tbox(2,1){0}
\tbox(3,1){1}
\tbox(4,1){1}
\tbox(5,1){2}
\tbox(6,1){2}
\tbox(7,1){2}
\tbox(8,1){2}
\tbox(9,1){2}

\tbox(0,2){0}
\tbox(1,2){0}
\tbox(2,2){0}
\tbox(3,2){0}
\tbox(4,2){0}
\tbox(5,2){0}
\tbox(6,2){1}
\tbox(7,2){1}
\tbox(8,2){2}
\tbox(9,2){2}

\tbox(0,3){0}
\tbox(1,3){0}
\tbox(2,3){0}
\tbox(3,3){0}
\tbox(4,3){0}
\tbox(5,3){0}
\tbox(6,3){0}
\tbox(7,3){0}
\tbox(8,3){0}
\tbox(9,3){2}

\tbox(0,4){1}
\tbox(1,4){1}
\tbox(2,4){1}
\tbox(3,4){1}
\tbox(4,4){0}
\tbox(5,4){1}
\tbox(6,4){1}
\tbox(7,4){1}
\tbox(8,4){1}
\tbox(9,4){2}

\tbox(0,5){0}
\tbox(1,5){1}
\tbox(2,5){1}
\tbox(3,5){0}
\tbox(4,5){0}
\tbox(5,5){1}
\tbox(6,5){0}
\tbox(7,5){1}
\tbox(8,5){0}
\tbox(9,5){1}

\endpspicture
\caption{Another representation of a permutation tableau} 
\label{PermTab2}
\end{figure}

Postnikov \cite{Postnikov1} has described a map that takes permutation
tableaux contained in a $k\times (n-k)$ rectangle to permutations in
$\cls_n$ with $k$ weak excedances.  In this paper we give a much simpler
description of this map, prove that it is a bijection, and show that
this map in fact preserves many more statistics.  Namely, the numbers
of 0's, 1's and 2's, respectively, in a permutation tableau equal
certain linear combinations of certain statistics defined on the
corresponding permutation.  Each of these statistics counts pairs of
indices $(i,j)$ in a permutation, according to the relative sizes of
the letters in those places and the place numbers themselves.  These
statistics were defined by Corteel in \cite{Corteel}.

We then define another bijection, taking permutations to permutations
and translating the statistics mentioned above into certain linear
combinations of generalized permutation patterns.  These combinations
between them contain precisely half the generalized patterns of length
3 with one dash (see Section~\ref{permpatts}).

We conclude our paper by giving various enumerative applications of
our bijections.  The structure of permutation tableaux is in many ways
more transparent than the structure of permutations, and hence lends
itself more easily to enumeration.  For example, by using our
bijections together with results of the second author \cite{Williams},
we are able to give the entire distribution of permutations according
to the number of occurrences of the generalized pattern $\bxca$.  This
is the first such result for any pattern of length 3 (or more).
However, this particular result, although first conjectured by the
present authors, was first proved by Corteel~\cite{Corteel}, whose
work provided us with a crucial piece of the puzzle solved here.

Another interesting consequence of the results presented here is that
the statistic counting permutation tableaux according to the number of
rows and number of 0's is an Euler-Mahonian statistic, that is, has
the same distribution as the bistatistic on permutations consisting of
the number of descents and the major index.  It follows that if 
we define $D_{k,n}(p,q,r)$ to be the polynomial enumerating the
permutation
tableaux contained in a $k \times (n-k)$ rectangle
 according to the number of $0$'s, $1$'s, and $2$'s, then
$D_{k,n}(p,1,1)$ is equal to Carlitz' classical $q$-analog $B_{n,k}(p)$
of the 
Eulerian numbers \cite{Carlitz}.  Additionally, 
$D_{k,n}(1,q,1)$ is equal to the more recent $q$-analog $E_{k,n}(q)$
of the 
Eulerian numbers that was studied in \cite{Williams}.


\section{Bijection from Permutation Tableaux to Permutations}
\label{sec-tab-perm}

In this section we describe a bijection $\Phi$ from permutation
tableaux to permutations.  More precisely, $\Phi$ is a bijection from
the set of permutation tableaux contained in a $k \times (n-k)$ rectangle
 to permutations in the symmetric group
$\cls_n$ with $k$ weak excedances.  Here, a {\em weak excedance} of a
permutation $\pi$ is a value $\pi(i)$ such that $\pi(i) \geq i$.  In
this situation we say that $i$ is a {\em weak excedance bottom} of
$\pi$ and that $\pi(i)$ is a {\em weak excedance top} of $\pi$.  To
make the notation less cumbersome, we abbreviate these as \emm
wexbottoms, and \emm wextops,, respectively.  {\em Non-weak excedance
bottoms} and {\em non-weak excedance tops} are defined in the obvious
way, in terms of $i$ and $\pi(i)$ such that $\pi(i)<i$, and are
abbreviated \emm non-wexbots, and \emm non-wextops,, respectively.
The \emm number, of weak excedances in $\pi$ will be denoted
$\wex\pi$.  Also, we let $\wexbotsum$ be the sum of all the wexbottoms
in $\pi$.

We remark that Postnikov \cite{Postnikov1} defined a map that is equal
to $\Phi$ but his description was more complicated and went
through the intermediate step of {\em web diagrams}.  Additionally,
his proof that it was a bijection used the geometry of the 
totally nonnegative part of the Grassmannian.  Another indirect
way to define the map $\Phi$ can be obtained from the Appendix of 
\cite{Williams2}: use the restriction from decorated permutations
to permutations of the composition
of the maps from Lemmas A.4 and A.5.

Our primary contributions in this section and the next are to: describe 
$\Phi$ in a simple enough way that we can actually give a 
combinatorial 
proofs of its properties;
prove that $\Phi$ translates a remarkable number of 
statistics on tableaux to statistics on permutations.

Before giving the bijection $\Phi$, we must define the {\em diagram}
$D(\T)$ associated with $\T$ as follows.  Regard the south-east border
of the partition $Y_\lambda$ contained in the $k \times (n-k)$
rectangle as giving a path (the {\it partition path}) 
$P = \{P_i\}_{i=1}^n$ of length $n$ from the northeast
corner of the rectangle to the southwest corner of the rectangle:
label each of the (unit) steps in this path with a number from $1$ to
$n$ according to the order in which the step was taken.  Then, remove
the $0$'s from $\T$ and replace each $1$ in $\T$ with a vertex.  We
will call the top vertex in each column a white vertex and all other
vertices black vertices.  Finally, from each vertex $v$, draw an
edge to the east and an edge to the south; each such edge should
connect $v$ to either a closest vertex in the same row or column, or
to one of the labels from $1$ to $n$.  The resulting picture is the
{\em diagram} $D(\T)$.  See Figure \ref{fig-diagram}.

We now define the permutation $\pi = \Phi(\T)$ via the following
procedure.  For each $i \in \{1, \dots , n\}$, find the corresponding
position on $D(\T)$ which is labeled by~$i$.  If the label $i$ is on a
vertical step of $P$, start from this position and travel straight
west as far as possible on edges of $D(\T)$. Then, take a ``zig-zag''
path southeast, by traveling on edges of $D(\T)$ south and east and
turning at each opportunity (i.e. at each new vertex).  This path will
terminate at some label $j\ge i $, and we let $\pi(i) = j$.  If $i$ is
not connected to any edge (equivalently, if there are no vertices in
the row of $i$) then we set $\pi(i)=i$.  Similarly, if the label $i$
is on a horizontal step of~$P$, start from this position and travel
north as far as possible on edges of $D(\T)$. Then, as before, take a
zig-zag path south-east, by traveling on edges of $D(\T)$ east and
south, and turning at each opportunity.  This path will terminate at
some label $j<i$, and we let $\pi(i) = j$.

See Figure \ref{zigzag} for a picture of the path taken by $i$.

\psset{unit=.7pt, arrowsize=7pt, linewidth=1pt}
\psset{linecolor=blue}
\newgray{grayish}{.90}
\newrgbcolor{embgreen}{0 .5 0}
\def\vblack(#1, #2)#3{\cnode*[linecolor=black](#1, #2){3}{#3}}
\def\vwhite(#1,#2)#3{\cnode[linecolor=black,fillcolor=white,fillstyle=solid](#1,#2){3}{#3}}
\countdef\x=23
\countdef\y=24
\countdef\z=25
\countdef\t=26

\def\tbox(#1,#2)#3{
\x=#1 \y=#2
\multiply\x by 36
\multiply\y by 36
\z=\x \t=\y
\advance\z by 36
\advance\t by 36
\psline(\x,\y)(\x,\t)(\z,\t)(\z,\y)(\x,\y)
\advance\x by 18
\advance\y by 18
\rput(\x,\y){{\bf #3}}}

\def\evbox(#1,#2)#3{
\x=#1 \y=#2
\multiply\x by 36
\multiply\y by 36
\advance\x by 10
\advance\y by 18
\rput(\x,\y){{\bf #3}}}

\def\ehbox(#1,#2)#3{
\x=#1 \y=#2
\multiply\x by 36
\multiply\y by 36
\advance\x by 18
\advance\y by -12
\rput(\x,\y){{\bf #3}}}

\def\contournumbers{
\evbox(3,0){5}
\evbox(4,1){3}
\evbox(4,2){2}
\evbox(4,3){1}
\ehbox(0,0){8}
\ehbox(1,0){7}
\ehbox(2,0){6}
\ehbox(3,1){4}
}

\def\bdot{\pscircle*{1.5mm}}
\def\wdot{\pscircle{1.5mm}}

\begin{figure}
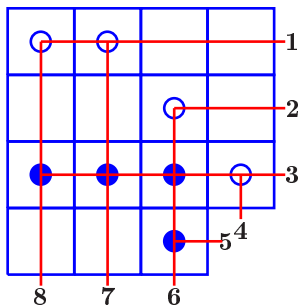

\pspicture(10,-10)(144,160)
\tbox(0,0){}
\tbox(1,0){}
\tbox(2,0){\bdot}

\tbox(0,1){\bdot}
\tbox(1,1){\bdot}
\tbox(2,1){\bdot}
\tbox(3,1){\wdot}

\tbox(0,2){}
\tbox(1,2){}
\tbox(2,2){\wdot}
\tbox(3,2){}

\tbox(0,3){\wdot}
\tbox(1,3){\wdot}
\tbox(2,3){}
\tbox(3,3){}

\contournumbers

\psset{linecolor=red}

\psline(150,126)(18,126)
\psline(18,126)(18,-6)

\psline(150,54)(18,54)
\psline(54,126)(54,-6)

\psline(90,90)(90,-6)
\psline(90,18)(116,18)

\psline(90,90)(150,90)
\psline(126,54)(126,30)

\endpspicture
\caption{\label{fig-diagram} The diagram of a tableau.  The topmost
$1$ in each column becomes a white vertex, and the other $1$'s become
black vertices.}
\end{figure}

\begin{figure}
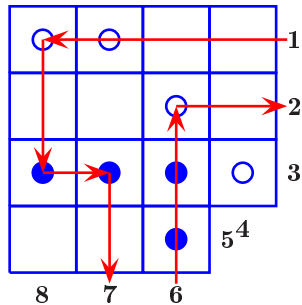

\pspicture(10,-10)(144,160)
\tbox(0,0){}
\tbox(1,0){}
\tbox(2,0){\bdot}

\tbox(0,1){\bdot}
\tbox(1,1){\bdot}
\tbox(2,1){\bdot}
\tbox(3,1){\wdot}

\tbox(0,2){}
\tbox(1,2){}
\tbox(2,2){\wdot}
\tbox(3,2){}

\tbox(0,3){\wdot}
\tbox(1,3){\wdot}
\tbox(2,3){}
\tbox(3,3){}

\contournumbers

\psset{linecolor=red}

\psline{->}(150,126)(18,126)
\psline{->}(18,126)(18,54)
\psline{->}(18,54)(54,54)
\psline{->}(54,54)(54,-6)

\psline{->}(90,-6)(90,90)
\psline{->}(90,90)(150,90)

\endpspicture
\caption{\label{zigzag}The paths taken by 1 and 6: $\pi(1)=7$, $\pi(6)=2$.}
\end{figure}

\begin{example}
If $\T$ is the permutation tableau whose diagram is
given in Figures \ref{fig-diagram} and \ref{zigzag}, then $\Phi(\T) =
74836215$.
\end{example}

The rest of this section will be devoted to proving various properties
of the map~$\Phi$, and in particular, that $\Phi$ is actually a bijection.

The following three lemmas are clear from the construction above.
\begin{lemma}\label{fp}
In $\Phi(\T)$, the letter $i$ is a fixed point if and only if there is
an entire row in $\T$ that has no 1's and whose right hand edge is
labeled by $i$.  In particular, $n,n-1,\ldots,n-m+1$ are fixed points
in $\pi$ if and only if the bottom $m$ rows of $\T$ (in the
$k\times(n-k)$ rectangle) consist entirely of 2's.
\end{lemma}

\begin{lemma}\label{E}
  Any directed step in a path on $D(\T)$ determines the path
  completely.
\end{lemma}

Lemma \ref{E} implies the following.

\begin{corollary}
$\Phi(\T)$ is a permutation.
\end{corollary}

\begin{lemma}\label{Lem3}
The weak excedance bottoms of $\pi = \Phi (\T)$ are precisely the
labels on the vertical edges of $P$.  The non-weak excedance
bottoms of $\pi$ are
precisely the labels on the horizontal edges of $P$.  In particular,
$\Phi(\T)$ is a permutation in $\cls_n$ with precisely $k$ weak
excedances.
\end{lemma}

We will now give some more definitions
and prove a refinement of Lemma \ref{fp}.

First we define a {\it relative permutation}: this is a {\it biword}
$\pi = \biw{a_1 \dots a_n},b_1\,\ldots\,b_n,$, where the $a_i$'s are
distinct integers and the $b_j$'s are also distinct. This denotes the
map sending $b_i \mapsto a_i$.  For example, the permutation 31524
corresponds to the biword $\biw 3\,1\,5\,2\,4, 1\,2\,3\,4\,5,$.  Note that it
is more common in the literature to write biwords in the other way,
that is, with the rows interchanged compared to our notation.  Because
we want to emphasize weak excedance tops and bottoms, we feel it is
more intuitive to have the tops in the top row and the bottoms in the
bottom row.  

We define the notion of {\it congruence} for relative permutations in
the obvious way, as follows.  Suppose that $a_{i_1} < a_{i_2} < \dots
< a_{i_n}$ and $b_{j_1} < b_{j_2} < \dots< b_{j_n}$; then define the
{\it reduction} of $\pi$ to be the permutation that one obtains by
replacing $a_{i_k}$ by $k$ and $b_{j_m}$ by $m$.  We now say that two
relative permutations are {\it congruent} if their reductions are
equal.

In keeping with the above definition, we define a {\it relative fixed point}
to be a
pair $b_p \mapsto a_p$ such that if $a_p$ is the $j$-th smallest letter
among the $a_i$ then $b_p$ is also the $j$-th smallest letter among
the $b_i$.

For example, in $\biw 5\,3\,1\,4,6\,2\,1\,3,$, the pair $2 \mapsto 3$ is a
relative fixed point, since each is the second smallest letter in its
row.

Note that we will use the biword notation $\biw{a_1 \dots
a_n},1\,\ldots\,n,$ as an alternative representation of the
permutation $a_1 \dots a_n$.

\begin{definition}\label{relative}
Fix a permutation tableau $\T$ and 
let $m:=n-k$ be the number of columns in $\T$.  We will 
construct a sequence of $m+1$ relative permutations
$\pi_0, \pi_1, \dots , \pi_m$ associated to $\T$ and will prove a 
generalization of Lemma \ref{fp} for these relative permutations.
Let $\pi_0 := \Phi(\T)$.  Let $r_1 < \dots < r_m$ be the 
non-excedence bottoms of $\pi_1$, in increasing order.
We construct $\pi_1$ by deleting the pair 
$r_1 \mapsto \pi_0(r_1)$ from the biword for $\pi_0$,
and then deleting any resulting relative fixed points.
In general, we construct $\pi_{i+1}$ from $\pi_i$ by 
deleting the pair $r_i \mapsto \pi_0(r_i)$ from the biword
for $\pi_i$, and then deleting any resulting relative fixed
points.
\end{definition}

\begin{example}\label{ex1}
Let $\T$ be the permutation tableau in Figure \ref{fig-diagram}.
Then we have that 
$\pi_0 =\biw 7\,4\,8\,3\,6\,2\,1\,5,1\,2\,3\,4\,5\,6\,7\,8,$,
$\pi_1 =\biw 7\,4\,8\,6\,2\,1\,5,1\,2\,3\,5\,6\,7\,8,$,
$\pi_2 =\biw 7\,8\,1\,5,1\,3\,7\,8,$,
$\pi_3 =\biw 7\,8\,5,1\,3\,8,$, 
and $\pi_4$ is the empty biword.
\end{example}

\begin{lemma}\label{interesting}
Use the notation of Definition \ref{relative}.
Fix $\T$ and choose some number $i$ which is less than the number
$m$ of columns of $\T$.  Construct a new permutation tableau
$\widetilde{\T}$
by deleting the rightmost $i$ columns of $\T$ and then deleting 
all of the resulting rows which contain only $0$'s.  Then 
$\Phi(\widetilde{\T})$ is congruent to $\pi_i$, as {\it relative permutations}.
\end{lemma}

\begin{example}
Let $\T$ be the permutation tableau in Figure \ref{fig-diagram}.
Then if we consider $\T$ and the resulting tableaux we get by 
cutting off the rightmost $i$ columns for $1 \leq i \leq 4$, the 
corresponding permutations that we get are:
$\biw 7\,4\,8\,3\,6\,2\,1\,5,1\,2\,3\,4\,5\,6\,7\,8,$,
$\biw 6\,3\,7\,5\,2\,1\,4,1\,2\,3\,4\,5\,6\,7,$,
$\biw 3\,4\,1\,2,1\,2\,3\,4,$,
$\biw 2\,3\,1,1\,2\,3,$, 
and the empty permutation.
Note that these are congruent to the relative permutations which 
we computed in Example \ref{ex1}.
\end{example}

\begin{proof}
By induction, it is enough to prove Lemma \ref{interesting} for 
the case $i=1$.  Suppose that $r$ is the smallest non-excedence
bottom of $\pi_0$, that is, $r$ is the number indexing the rightmost
column of $\T$, and suppose that $r \mapsto h$ in $\pi_0$
(for some $h < r$).  Then by definition, $\pi_1$ is equal to the 
result of deleting $r \mapsto h$ and any resulting relative 
fixed points from the relative permutation $\pi_0$.

Let $\widetilde{\T}$ be the permutation tableau formed by deleting
the rightmost column of $\T$ and then deleting any resulting
all-zero rows.  Since it will not affect the congruence class of 
$\Phi(\widetilde{\T})$, we can label the partition path of 
$\widetilde{\T}$ with the labels inherited from $\T$; in other words,
we can label the new partition path with the numbers 
$\{1, \dots , n\} \setminus \{r, j_1, j_2, \dots \}$, where the 
$j_i$'s index the all-zero rows that we deleted.  Now let us
analyze the difference between $\Phi(\widetilde{\T})$ and $\pi_0$.
By consideration of the map $\Phi$, $\Phi(\widetilde{\T})$ is 
identical to $\pi_0$ except that:
\begin{itemize}
\item we have removed $r \mapsto h$ from the biword for $\pi_0$ 
\item we have removed resulting relative fixed points 
\item if in $\pi_0$ we had $i \mapsto j$ where either $j=r$ or 
  $j$ indexed a row that is not present in $\widetilde{\T}$, then
in $\Phi(\widetilde{\T})$ we have that $i \mapsto g$, where 
$g$ is the maximal label in 
$\{1, \dots , n\} \setminus \{r, j_1, j_2, \dots \}$ which is 
less than $j$.  
\end{itemize}
It is clear now that $\Phi(\widetilde{\T})$ is congruent to $\pi_1$.
\end{proof}

\begin{lemma}\label{Rfp}
Use the notation of Definition \ref{relative}.
After we delete   
$r_i \mapsto \pi_0(r_i)$ from the biword for $\pi_{i-1}$, 
we will have a relative fixed point
$j \mapsto  j'$ in the resulting relative permutation if and only if
every entry in row $j$ of $\T$ to the left of the column indexed by 
$r_i$ is a zero.
\end{lemma}

\begin{proof}
We will prove this lemma by induction.  Suppose that it is true for 
$i \leq I$; we will now prove it for $I+1$.
Let $\widetilde{\T}$ be the tableau obtained by deleting the rightmost
$I$ columns of $\T$ and then deleting any resulting all-zero rows.
By Lemma \ref{interesting}, we have that 
$\Phi(\widetilde{\T})$ is congruent to $\pi_I$, and by the induction
hypothesis, these relative permutations have no relative fixed points.

Suppose that with the exception of the rightmost column of 
$\widetilde{\T}$ (which corresponds to the column in $\T$ indexed by
$r_{I+1}$), all entries in row $j$ of $\widetilde{\T}$ are $0$.
Label the partition path of $\widetilde{\T}$ with the consecutive numbers
$1, \dots , n'$ for some $n'$ and suppose that $r$ is the label indexing
the  
rightmost column.  Clearly $r > j$.  Since $\pi_I$ has no relative fixed
points, the entry of $\widetilde{\T}$ in column $r$ and row $j$
must be a $1$, and additionally every entry in the same column below
this $1$ must also be a $1$.  It follows from the definition of 
$\Phi$ that in the permutation 
$\Phi(\widetilde{\T})$, we have $j \mapsto j+1$.  Again by the 
definition of $\Phi$, we have $r \mapsto h$ in 
$\Phi(\widetilde{\T})$ for some $h < j$.  Recalling that 
$\Phi(\widetilde{\T})$ is congruent to $\pi_I$ and $r > j$, 
it follows that when we delete $r \mapsto h$ in $\pi_I$, 
$j \mapsto j+1$ will become a relative fixed point in the resulting
relative permutation.  

All of these steps can be reversed, so we are done.
\end{proof}

We are finally ready to prove the following theorem.

\begin{theorem}\label{thm-phi}
  The map $\Phi$ is a bijection from permutation
  tableaux to permutations.
\end{theorem}

\begin{proof}
To prove that $\Phi$ is a bijection, we will give an explicit
description of its inverse, by reverse-engineering $\Phi$
so as to be consistent with Lemma
\ref{Rfp}.

Let $\pi$ be the permutation
$\biw{a_1 \dots a_n},1\,\ldots\,n,$.
We now give the procedure for computing
$\Phi^{-1}(\pi)$, column by column, from right to left.

\begin{itemize}
\item[0.] 
 Compute the weak excedance bottoms of $\pi$ to get the shape of the
 partition in $\T$ (see Lemma~\ref{Lem3}).
 Let $\tilde{\pi}=\pi$.
  
\item[1.] Check for relative fixed points in $\tilde{\pi}$.  If
$j \mapsto j'$ is a relative fixed point then by Lemma \ref{Rfp}, 
$\T$ must have 0's in the as-yet-undetermined 
part of the row corresponding to the weak excedance
bottom $j$; fill in these entries.
Recompute $\tilde{\pi}$ by removing the relative fixed
points.
  
\item[2.] Suppose we have determined the content of the $i$ rightmost
  columns (but have not determined the content of the other columns).
  Then look at the next column to the left, which is indexed
  by a non-excedance bottom $r$ (that is, by the label on the
  horizontal step at the bottom of that column (Lemma~\ref{Lem3})).
  Knowing that $r \to a_r$ in $\pi$ uniquely determines the position
  $p$ of the highest $1$ in the column corresponding to $r$, since
  there is a unique zig-zag path going backwards (north-west) from
  $a_r$ to a box in the column above~$r$.  Insert a $1$ at that
  position and $0$'s in all boxes above it which are in the same
  column.  Also, insert $1$'s into all undetermined boxes below $p$.
  (Note that we know that all nonzero boxes below position $p$ must
  also be $1$'s; otherwise, if there were some $0$ below the
  $1$ then everything to its left would have to be a $0$,
  and then by Lemma \ref{Rfp}, we would have a relative fixed point
  in $\tilde{\pi}$, contradiction.)
  Reduce $\tilde{\pi}$ by removing the column $r \mapsto a_r$ from the
  biword for~$\tilde\pi$.  Go to step 1.
\end{itemize}

It is now clear from this construction that our resulting tableau $\T$
will be a permutation tableau, and moreover, that it will be the 
inverse image $\Phi^{-1}(\pi)$.
\end{proof}

\begin{example}
  Let $\pi = 514263$.  Since $\pi$ is in $S_6$ and has three weak
  excedances and $a_6\ne6$, our permutation tableau $\T$ will be
  contained in a $3 \times 3$ rectangle (the resulting tableau is
  shown in Figure~\ref{inv-tableau}).  As in step 0, we want to first
  compute the shape of the associated partition.  Since $1, 3, 5$ are
  the wexbottoms, and $2,4,6$ are the non-wexbottoms, this uniquely
  determines a path (the {\it partition path}) from the northeast
  corner of the rectangle to the southwest corner of the rectangle
  with vertical steps in positions $1,3,5$ and horizontal steps in
  positions $2,4,6$.  That is, our partition has the shape $(3,2,1)$.
  We now draw this partition, labeling the edges of its southeast
  border accordingly with the numbers $1, \dots, 6$, and set
  $\tilde{\pi} = (5,1,4,2,6,3)$.
  
  Going to step 1, we see that $\tilde{\pi}$ has no relative fixed
  points.
  
  Going to step 2, the fact that $2 \to 1$ in $\pi$ implies that the
  rightmost column (which consists of a single box) contains a 1
  in the top row.  We now reduce the permutation $\tilde{\pi} =
  \biw5\,1\,4\,2\,6\,3,1\,2\,3\,4\,5\,6,$ by removing $(2,1)$,
  obtaining $\tilde{\pi} = \biw5\,4\,2\,6\,3,1\,3\,4\,5\,6,$.
  
  Going back to step 1, we see that there are no relative fixed points
  in $\tilde{\pi}$.
  
  Going to step 2, since $4 \to 2$ in $\tilde\pi$ it is clear that
  the highest box in the column indexed by $4$ must contain a 1.
  All undetermined boxes below this 1 must contain 1's also.
  We now reduce the permutation $\tilde{\pi} =
  \biw5\,4\,2\,6\,3,1\,3\,4\,5\,6,$ by removing $4 \to 2$, obtaining
  $\tilde{\pi} = \biw5\,4\,6\,3,1\,3\,5\,6,$.
  
  Going back to step 1, we now see that $\tilde{\pi}$ has the relative
  fixed point $(3,4)$.  Therefore the undetermined part of the row
  corresponding to $3$ consists of zeros.  Now we reduce the
  permutation $\tilde{\pi} = \biw5\,4\,6\,3,1\,3\,5\,6,$, obtaining
  $\tilde{\pi} = \biw5\,6\,3,1\,5\,6,$.
  
  Going to step 2, since $6 \to 3$ in $\tilde\pi$ the top box in the
  column corresponding to $6$ has a 1.  All undetermined boxes
  below that contain 1's.  We have now filled in all columns of 
 the tableau---obtaining the permutation tableau
 in Figure~\ref{inv-tableau}---so we are done.

\end{example}

\psset{unit=.8pt, linewidth=1pt}

\begin{figure}
\def\evbox(#1,#2)#3{
\x=#1 \y=#2
\multiply\x by 36
\multiply\y by 36
\advance\x by 7
\advance\y by 18
\rput(\x,\y){{\bf #3}}}

\def\ehbox(#1,#2)#3{
\x=#1 \y=#2
\multiply\x by 36
\multiply\y by 36
\advance\x by 18
\advance\y by -7
\rput(\x,\y){{\bf #3}}}

\pspicture(10,-10)(144,160)
\tbox(0,0){1}

\tbox(0,1){0}
\tbox(1,1){1}

\tbox(0,2){1}
\tbox(1,2){1}
\tbox(2,2){1}

\evbox(1,0){5}
\evbox(2,1){3}
\evbox(3,2){1}
\ehbox(0,0){6}
\ehbox(1,1){4}
\ehbox(2,2){2}

\endpspicture
\caption{\label{inv-tableau}The permutation tableau for $\pi=514263$}
\end{figure}


\section{How $\Phi$ translates statistics}

The six permutation statistics in the following definition will be
related to the statistics recording the numbers of 0's, 1's and 2's in
permutation tableaux.  The first four of these refine Postnikov's
definition of alignment \cite{Postnikov} (see \cite{Williams}); all of
these statistics were defined by Corteel~\cite{Corteel}.

\begin{figure}[h]
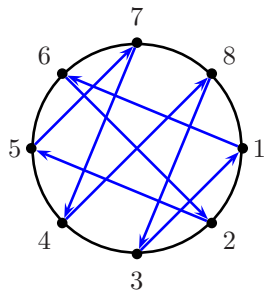

\begin{center}
\psset{unit=1pt, arrowsize=4pt, linewidth=1pt}
\pspicture(-60, -60)(60,60)
\pscircle[linecolor=black](0,0){40}
\cnode*[linewidth=0, linecolor=black](40,0){2}{1}
\cnode*[linewidth=0, linecolor=black](28.28,-28.28){2}{2}
\cnode*[linewidth=0, linecolor=black](0, -40){2}{3}
\cnode*[linewidth=0, linecolor=black](-28.28,-28.28){2}{4}
\cnode*[linewidth=0, linecolor=black](-40,0){2}{5}
\cnode*[linewidth=0, linecolor=black](-28.28,28.28){2}{6}
\cnode*[linewidth=0, linecolor=black](0,40){2}{7}
\cnode*[linewidth=0, linecolor=black](28.28,28.28){2}{8}
\rput(46.30,0){$1$}
\rput(35,-35){$2$}
\rput(0,-50){$3$}
\rput(-35,-35 ){$4$}
\rput(-46.30,0){$5$}
\rput(-35,35){$6$}
\rput(0, 50 ){$7$}
\rput(35,35){$8$}
\ncline{->}{1}{6}
\ncline{->}{2}{5}
\ncline{->}{3}{1}
\ncline{->}{4}{8}
\ncline{->}{5}{7}
\ncline{->}{6}{2}
\ncline{->}{7}{4}
\ncline{->}{8}{3}

\endpspicture
\end{center}
\caption{A chord diagram for the permutation $65187243$} 
\label{chorddiagram}
\end{figure}

\begin{definition}\label{def-ac}
Given a permutation $\pi=a_1a_2\ldots a_n$, let
\begin{eqnarray*}
\aa(i) &=& \{j \st j<i\le a_i < a_j \},\\
\ab(i) &=& \{j \st a_j < a_i < i < j \},\\
\ac(i) &=& \{j \st j\le a_j < a_i < i \},\\
\ad(i) &=& \{j \st a_i< i < j\le a_j \},\\
\ca(i) &=& \{j \st j<i\le a_j < a_i \},\\
\cb(i) &=& \{j \st a_i < a_j < i < j \}.
\end{eqnarray*}
We then set
\begin{eqnarray*}
\aa(\pi) = \sum_i{|\aa(i)|},
\end{eqnarray*}
and likewise for the other five statistics.
\end{definition}
Observe that if we draw the permutation as a chord diagram on a
circle, as in Figure \ref{chorddiagram}, then $j \in A_{**}(i)$ means
that the chords starting at $i$ and $j$ do not intersect and roughly
``point in the same direction'' (see \cite{Williams} for more details); 
we will say that this is an {\em
alignment of type $\mathrm{A_{**}}$}.  And if $j
\in\mathrm{C_{**}}(i)$ then the chords starting at $i$ and $j$ cross
each other; we will say that this is a {\em crossing of type
$\mathrm{C_{**}}$}.  Note that the subscripts in our notation refer to
whether the positions $i$ and $j$ are wexbottoms or non-wexbottoms of
the permutation.  For example, in Figure \ref{chorddiagram}, the
chords beginning at $3$ and $5$ form an alignment of type $\Ad$, and
the chords beginning at $2$ and $4$ form a crossing of type~$\Cb$.

\begin{theorem}\label{equidistribution1}
  Let $T(k,a,b,c)$ be the set of permutation tableaux with $k$ rows,
  $(n-k)$ columns, $a$ 0's, $b$ 1's and $c$ 2's.  Let $M(k,a,b,c)$ be
  the set of all permutations $\pi \in\cls_n$ with
\begin{itemize}
\item%
 $k = \wex(\pi)$,
\item%
  $a =                  
  \Aa(\pi)+\Ab(\pi)+\Ac(\pi)$,
\item%
  $b = 
 \Ca(\pi)+\Cb(\pi)+(n-k)$,
\item%
  $c =\Ad(\pi)$.
\end{itemize}
Then $|T(k,a,b,c)| = |M(k,a,b,c)|$.  Moreover, the 
map $\Phi$ is a bijection from $T(k,a,b,c)$ to $M(k,a,b,c)$
such that the weak excedence bottoms of $\pi = \Phi(\T)$
are precisely the labels on the vertical edges of the partition
path $P$ associated with $\T$.
\end{theorem}

Because of Lemma \ref{Lem3}, it is enough to prove
that $\Phi$ is a bijection.
This will be done using a lemma of Corteel \cite{Corteel},
and 
Propositions \ref{twos-Ad}
and \ref{black} below. 

\begin{lemma}[Corteel \cite{Corteel}]\label{lemma-corteel}
Let $k$, $n$, $\aa(\pi)$, $\ab(\pi)$, $\ac(\pi)$, 
$\ad(\pi)$, $\ca(\pi)$, $\cb(\pi)$ be as above.
Then
  $$
  \aa(\pi)+\ab(\pi)+\ac(\pi)+\ad(\pi)+\ca(\pi)+\cb(\pi) = (k-1)(n-k).
  $$
\end{lemma}

\begin{prop}\label{twos-Ad}
If $\Phi(\T) =\pi$ then the number of $2$'s in $\T$ is equal to the
number of 
alignments of type $\Ad$  in $\pi$.
\end{prop}

\begin{proof}
Recall that if $\pi=\Phi(\T)$ then the wexbottoms and the
   non-wexbottoms of~$\pi$ correspond to the labels of the vertical
   and horizontal steps, respectively, in the south east border of the
   partition underlying $\T$.  Note that the position of every 2 in
   $\T$ can be given by specifying the label of the edge above it and
   the edge to its left.  The label $i$ of the edge above it will be a
   non-wexbottom, and the label $j > i$ of the edge to its left will
   be a weak excedance bottom.  Since $j > i$, and $j$ is a wexbottom,
   and $i$ is a non-wexbottom, the pair $(i,j)$ is precisely an
   alignment of type $\Ad$.  Conversely, any alignment of type $\Ad$
   is a pair $(i,j)$ where $i<j$, and $i$ is a non-wexbottom, and $j$
   is a wexbottom.  This implies that $i$ is the label of a horizontal
   step, and $j$ is the label of a vertical step.  The fact that $i<j$
   implies that the box of the tableau indexed as above by $i$ and $j$
   contains a 2.
\end{proof}

\begin{prop}\label{black}
  Under the bijection $\Phi$, there is a one-to-one correspondence
  between black vertices in the diagram of the permutation tableau,
  and crossings of types $\Ca$ and $\Cb$ in the permutation.
\end{prop}

Before proving this proposition, we will illustrate the main idea
with an example.

\begin{example}
Again consider the tableau $\T$ in Figure \ref{fig-diagram}.
This tableau corresponds to the permutation 
$\pi = 74836215$, which has a total of four crossings.
The pair of chords $2 \mapsto 4$, $3 \mapsto 8$, and the pair of 
chords $1 \mapsto 7$, $3 \mapsto 8$ are crossings of type
$\ca$, while  the pair of chords
$6 \mapsto 2$, $8 \mapsto 5$ and the pair of chords
$7 \mapsto 1$, $8 \mapsto 5$ are crossings of type $\cb$.

Observe that under the bijection $\Phi$, the
paths $(2\to 4)$ and $(3 \to 8)$ intersect in a unique horizontal
edge: the edge between the horizontally adjacent black 
and white vertices in row $3$ of the tableau.  We will associate
this crossing to the black vertex which is the leftmost of the 
two vertices.

Similarly, the paths $(1\to 7)$ and $(3 \to 8)$ intersect 
in a unique horizontal edge: the edge between the two leftmost black
vertices in row $3$.  We will associate this crossing to the black
vertex which is the leftmost of the two vertices.

On the other hand, the paths $(6 \to 2)$ and $(8 \to 5)$ intersect
in the vertical edge between the two black vertices in the column indexed by 
$6$.  We will associate this crossing to the black vertex which is 
the bottom of these two vertices.

Simiarly, the paths $(7 \to 1)$ and $(8 \to 5)$ intersect in the 
vertical edge between the two vertices in the column indexed by $7$.
We will associate this crossing to the black vertex which is the 
bottom of these two vertices.

In this way, the crossings of $\pi$ get associated bijectively
with the black vertices of $\T$.
\end{example} 

We will now make rigorous the idea that the above example suggests.

\begin{proof}
  Recall that black vertices correspond to those 1's in a tableau that
  are not topmost in their columns.  Let $D$ be the diagram of $\T$,
  and let $\pi=(a_1, \dots , a_n)$ be $\Phi(\T)$.  We will construct a
  map $\phi$ (induced by $\Phi$) which takes each crossing $(i,j)$
  (where $i<j$) of type $\Ca$ or $\Cb$ in $\pi$ to a black vertex $d$
  in $D$, and show that this is a bijection.  The map $\phi$ is
  defined as follows.  Let $(i,j)$ be a crossing of type $\Ca$ or
  $\Cb$.  We claim that the paths $(i\to a_i)$ and $(j \to a_j)$
  intersect in a unique edge.  If that edge is horizontal, then let
  $d$ be the left vertex of the edge.  If that edge is vertical, then
  let $d$ be the bottom vertex of the edge.
  
  First we need to show that the paths $(i\to a_i)$ and $(j \to a_j)$
  intersect in an edge.  We will prove this when $(i,j)$ is a crossing
  of type $\Ca$; the proof for $\Cb$ is similar.  Since $i<j$ and $a_i
  < a_j$, it is clear that the paths must cross each other at least
  once.  
  
  Consider the first point $x$ at which the path $(j \to a_j)$
  intersects the path $(i \to a_i)$.  We will show that the
  intersection here will contain an edge.  Clearly this intersection
  must be in the zig-zag portion of the path $(i \to a_i)$. If we let
  $d_1, d_2, \dots , d_t$ be the sequence of vertices encountered by
  the path $(i \to a_i)$ in its zig-zag portion, then, by construction
  of that path, there are no vertices in the diagram $D$ between any
  $d_r$ and $d_{r+1}$.  Note that if the path $(j \to a_j)$ intersects
  the path $(i \to a_i)$ in only the point $x$ (rather than an edge
  containing $x$), then it is easy to see---using condition (2) in the
  definition of permutation tableaux---that $x$ must actually be a
  vertex in $D$, located between some $d_r$ and $d_{r+1}$.  This is a
  contradiction.
  
  Next, we show that the paths $(i\to a_i)$ and $(j \to a_j)$
  intersect in a {\em unique} edge.  If the two paths were to
  intersect a second time (and they may indeed intersect again in a
  {\em vertex}), then this intersection must take place in the zig-zag
  portion of {\em both} paths.  Such a point $e$ of intersection must
  be approached via a south step by $(i \to a_i)$ and must be
  approached via an east step by $(j \to a_j)$.  But then, according
  to the procedure defining $\Phi$, the path $(i \to a_i)$ will
  immediately turn east, and the path $(j \to a_j)$ will immediately
  turn south. Therefore this intersection is not an edge intersection.
  
  We have thus shown that $\Phi$ induces a well-defined map from
  crossings to black vertices.  We will now show that this map is a
  bijection by constructing its inverse.  Namely, to each black vertex
  in $D$ we need to produce a crossing of type $\Ca$ or $\Cb$.  We do
  this as follows.  Given a black vertex $d$, there is a path $(i\to
  a_i)$ on~$D$ which enters $d$ by going south, and then leaves $d$
  going east.  (It is easy to see that such a path exists by tracing
  backwards through the algorithm that defined the map $\Phi$.)
  
  If the path $(i\to a_i)$ is an excedance, then consider the unique
  path $(j \to a_j)$ which enters $d$ traveling west.  This path must
  be a weak excedance, as it is only the paths of weak excedances
  which contain steps to the west.  Moreover, $(i,j)$ must form a
  crossing of type $\Ca$, since the two paths intersect in an edge
  (and we have seen that two paths which are both weak excedances may
  not intersect in an edge more than once).
  
  On the other hand, if the path $(i\to a_i)$ is a non-excedance, then
  consider the unique path $(j \to a_j)$ which enters $d$ traveling
  north.  Clearly this path must be a non-excedance, as it is only the
  paths of non-excedances which contain steps north.  Moreover, $(i,j)$
  must form a crossing of type $\Cb$, since the two paths must
  intersect in a unique edge.
  
  Therefore $\phi$ is a bijection between the set of $\Ca$- and
  $\Cb$-crossings in $\pi$, and the set of black vertices in $D$.
\end{proof}

We now finish the proof of Theorem \ref{equidistribution1} with the 
following argument.

\begin{proof}
The first, third, and fourth parts of the theorem follow from 
Lemma \ref{Lem3}, Proposition \ref{twos-Ad}, and Proposition \ref{black}, 
respectively.  It remains to prove the second part of the theorem.
Let $m=\Aa(\pi)+\Ab(\pi)+\Ac(\pi)$.  We know that $a+b+c=k(n-k)$.  
By Lemma \ref{lemma-corteel},
we have that $\Aa(\pi)+ \Ab(\pi)+\Ac(\pi)+\Ad (\pi)+
\Ca(\pi)+\Cb(\pi) = (k-1)(n-k)$.  Therefore $m+c+b-(n-k) =
(k-1)(n-k)$, which implies that $m+c+b = k(n-k) = a+b+c$, and hence
$m=a$.
\end{proof}

\section{Permutation patterns}\label{permpatts}

In this section we introduce necessary terminology and definitions
that will be used in the next section, where we construct a bijection
$\Psi:\cls_n\ra\cls_n$.  This bijection proves the equidistribution of
certain linear combinations of the statistics in
Definition~\ref{def-ac} (alignments and crossings) with certain linear
combinations of generalized permutation patterns, which we define
below.  The composition of $\Psi$ and the bijection $\Phi$ from
Section~\ref{sec-tab-perm} then proves the equidistribution of our
tableaux statistics (numbers of 0's, 1's and 2's) with the pattern
statistics to be defined here.

A \emm classical permutation pattern, $p=p_1p_2\ldots p_k$ is simply a
permutation, and an \emm occurrence, of $p$ in a permutation
$\pi=a_1a_2\ldots a_n$ is a subsequence
$a_{i_1},a_{i_2},\ldots,a_{i_k}$ of $\pi$ (where $i_1<i_2<\cdots<i_k$)
whose letters are in the same relative order as in $p$.  For example,
the permutation 416235 has two occurrences of the pattern $\bxcxa$,
namely the subsequences 462 and 463.

In the literature, the pattern $\bxcxa$ is usually denoted simply by
$231$.  We write it here with dashes between consecutive letters in
order to emphasize that there are no restrictions on the distance
between the letters in a permutation that form an occurrence of the
pattern.  A \emm generalized pattern, is a pattern where some pairs of
adjacent letters may lack a dash between them.  Such an absence
indicates that the corresponding letters must be adjacent in an
occurrence of the pattern in a permutation.  For example, the pattern
$\bxca$ occurs only once in 416235, namely as 462.  In the subsequence
463, whose letters are in the same relative order as those of $\bxca$,
the last two letters are not adjacent in 416235 as required for an
occurrence of $\bxca$.  

Generalized patterns were first introduced systematically by Babson
and Steingrimsson in \cite{bast}, but some instances had been treated
previously in various contexts.  For example, the pattern $\caxb$ is
implicit in \cite{foata-zeil} and in \cite{simstan} (where the similar
patterns $\bxca$, $\acxb$ and $\bxac$ are also treated), and dashless
patterns, such as $\abc$, appeared already in \cite{goulden-jackson}.

The reason for writing patterns in parentheses is that we will
consider them as \emm functions, from the set of permutations to the
natural numbers, where the value of a pattern $p$ on a permutation
$\pi$ is the number of different occurrences of $p$ in $\pi$.  For
example, if $\pi=416235$, as above, then $\bxcxa\pi=2$ and
$\bxca\pi=1$.

It is easy to see that there are exactly twelve different patterns of
length 3 with one dash.  Six of these will be considered here, namely
$\axcb$, $\bxca$, $\cxba$, $\baxc$, $\caxb$ and $\cbxa$.  These are
all the patterns of length 3 with one dash whose two letters \emm
$not$, separated by a dash are in decreasing order.

A \emm descent, in a permutation $\pi=a_1a_2\ldots a_n$ is an $i$ such
that $a_i>a_{i+1}$. We say that $a_i$ is a \emm descent top, and
$a_{i+1}$ a \emm descent bottom,.  The set of descent tops is denoted
$\destopset$, and the set of descent bottoms $\desbotset$.  Moreover,
we let $\destopsum$ be the sum of the elements of $\destopset$, and
likewise for $\desbotsum$.

We now define the linear combinations of patterns whose joint
distribution on permutations matches the distribution of 0's, 1's and
2's on permutation tableaux.

\newcommand\aaa{a} \newcommand\bbb{b} \newcommand\ccc{c}
\begin{definition}
  Given a permutation $\pi$, let
\begin{eqnarray*}
\aaa(\pi) &=& \baxc\pi + \cxba\pi + \caxb\pi -\ch \des\pi,2,,\\
\bbb(\pi) &=& \bxca\pi + n - 1 - \des\pi,\\
\ccc(\pi) &=& \axcb\pi + \cbxa\pi -\ch \des\pi,2,.
\end{eqnarray*}
\end{definition}

It is important to note that since we will be considering the
%
%
quadruple statistic consisting of $\aaa$, $\bbb$, $\ccc$ and the
number of descents, the terms $\ch\des\pi,2,$ and $n-1-\des\pi$ in the
above definition only effect a \emm shift, of the statistics involved,
but not an essential modification.  For example, the bistatistic
$(\des,\bbb)$ can be seen as a 2-dimensional array of numbers, with
the $k$-th entry in the $i$-th row consisting of the number of
permutations $\pi$ with $i$ descents and with $\bbb(\pi)=k$.  Thus,
replacing $\bbb=\bxca + n - 1 - \des\pi$ ~by~ $\bbb'=\bxca$, we would
only shift the nonzero entries in each row $n-1-i$ steps to the left.
A similar, but more complicated, statement is true for the quadruple
statistic $(\des,\aaa,\bbb,\ccc)$.  The upshot of this is that if we
delete the terms $\ch\des\pi,2,$ and $(n-1-\des\pi)$ from the
definitions of $\aaa$, $\bbb$, and $\ccc$, the resulting quadruple
statistic changes only in that we are disregarding some ``initial''
zeros.

It is also important to note that the sum $\axcb\pi+\cbxa\pi$ in
$\ccc$ is equal to $\desbotsum\pi-\des\pi$.  Namely, for each descent
$\ldots yx\ldots$ in $\pi$, the pattern $\axcb$ counts the letters to
the left of the descent that are smaller than its descent bottom
($x$).  The letters to the right of the descent, and smaller than $x$,
are counted by $\cbxa$, so clearly we are counting \emm all,\/ letters
in $\pi$ that are smaller than~$x$.  Analogously, the sum
$\baxc\pi+\cxba\pi$ in $\aaa$ equals the sum of $n-t$ over all descent
tops $t$ in $\pi$.

This leaves $\caxb$ in $\aaa$, which sums the \emm left embracing
numbers, in $\pi$, so called because $\caxb$ counts, for each letter
$x$ in $\pi$, the descents to the left of $x$ that embrace $x$, that
is, where the letters of the descent are larger and smaller,
respectively, than $x$.  Analogously, the pattern $\bxca$ in $\bbb$
sums the \emm right embracing numbers, in~$\pi$.

To be more precise, we define the right embracing number of each
letter $\ell$ in $\pi$, denoted $\rembr(\ell)$, as the number of
descents $\ldots yx\ldots$ to the right of $\ell$ in $\pi$ such that
$x<\ell<y$.

\begin{lemma}\label{patt-const}
  Let $\des$ be the number of descents in a permutation $\pi$, and let
  $a(\pi)$, $b(\pi)$ and $c(\pi)$ be as above.  Then
  $$
  \aaa(\pi)+\bbb(\pi)+\ccc(\pi) = (\des+1)(n-\des-1).
  $$
\end{lemma}
\begin{proof}
  Each of the patterns involved in $\aaa+\bbb+\ccc$ counts certain
  letters to the left or to the right of each descent in $\pi$.
  Together they count, for each descent in $\pi$, all the letters in
  $\pi$ not belonging to the descent itself.  There are, of course,
  $n-2$ such letters for each descent.  Thus, the sum of all the
  patterns in $\aaa(\pi)+\bbb(\pi)+\ccc(\pi)$ is $\des\cdot(n-2)$.
  Completing the proof now only requires a routine calculation.
\end{proof}

\section{Another bijection}

We now describe the construction of a bijection $\Psi: \cls_n \ra
\cls_n$ that takes a permutation $\pi$ to a permutation $\tau$ such
that the set of descent tops in $\pi$ determines the set of weak
excedance tops in $\tau$ and the set of descent bottoms in $\pi$
determines the set of weak excedance bottoms in $\tau$. Moreover, the
right embracing number of $i$~in~$\pi$ becomes $\ca(i)$ in $\tau$ if
$i$ is a wexbottom in $\tau$ and becomes $\cb(i)$ in $\tau$ otherwise.

This bijection is based on the same idea as the central bijection
in~\cite{csz}, which in turn was shown, in \cite{csz}, to be
essentially equivalent to several seemingly different bijections in
the literature, due to Foata--Zeilberger \cite{foata-zeil}, Fran\c
con--Viennot \cite{fravie}, de~M\'edicis--Viennot \cite{medvie} and
Simion--Stanton \cite{simstan}, respectively.  More precisely, our
bijection $\Psi$ here uses the same \emm data, as the
bijection $\Phi$ in \cite{csz}, and these data (descent tops and
bottoms and right embracing numbers) completely determine a
permutation.  However, our bijection $\Psi$ uses this data in a
different way than $\Phi$ in \cite{csz}, and produces a different
permutation, turning descent tops into weak excedance tops etc.,
whereas $\Phi$ in \cite{csz} turns descent tops into excedance tops,
etc.

Recall the biword notation for permutations.  
For example, we write the permutation 31524 as
\begin{equation*}
\biw3\,1\,5\,2\,4,1\,2\,3\,4\,5,.
\end{equation*}
In order to construct $\Psi (\pi)$ (where $\pi \in \cls_n$),
we first construct two biwords, $\biw
f',f,$ and $\biw g',g,$, and then form the biword $\tau'=\biw
f'\;g',f\;g,$ by concatenating $f$ and $g$, and $f'$ and $g'$,
respectively.  The words $f,f',g,g'$ are defined as follows (we will
prove later, in Theorem~\ref{phi-thm}, that this is possible):
\begin{itemize}
\item%
  The letters of $f$ consist of the set obtained by adding 1 to each
  of the descent bottoms in $\pi$ and then adjoining the letter 1.
  The letters of $f$ are ordered increasingly.  {\it These 
  letters will be the wexbottoms of $\tau$.}

\item%
  The letters of $g$ consist of the set obtained from the non-descent
  bottoms in $\pi$ by removing the letter $n$ and adding 1 to the
  remaining letters.  The letters of $g$ are ordered increasingly.
  {\it These letters will be the non-wexbots of $\tau$.} 

\item%
  The letters of $f'$ consist of the set obtained by subtracting 1
  from each of the descent tops of $\pi$ and then adjoining the letter
  $n$.  The letters of $f'$ are ordered so that, for $a$ in $f'$,
  $\ca(a)$ in $\tau$ is the right embracing number of $a$ in~$\pi$.
  (Observe that $\ca(a)$ only depends on the relative order of the
  wextops in $\tau$, together with their corresponding wexbottoms.)
 {\it These letters will be the wextops of $\tau$.}

\item%
  The letters of $g'$ consist of the set obtained by removing $1$ from
  the set of non-descent tops in $\pi$ and then subtracting 1 from the
  remaining letters.  The letters of $g'$ are ordered so that, for $a$
  in $g'$, $\cb(a)$ is the right embracing number of $a$ in~$\pi$.
  {\it These letters will be the non-wextops of $\tau$.}
\end{itemize}

Rearranging the columns of $\tau'$, so that the bottom row is in
increasing order, we obtain the permutation $\tau=\Psi(\pi)$ as the
top row of the rearranged biword.  Before we prove that this can
always be done in the way described, we give an example.

Let $\pi=215896374$.  Then $\pi$ has
\smallskip
\begin{verbatim}
    Descent bottoms:  1 3 4 6     Non-descent bottoms:  2 5 7 8 9

    Descent tops:     2 6 7 9     Non-Descent tops:     1 3 4 5 8
\end{verbatim}
\smallskip
The right embracing numbers are 2 for 5, 1 for 6 and 8, and 0 for all
others:
\bigskip
{\setlength\baselineskip{.7\baselineskip}
\begin{verbatim}
                          21-5-8-963-74
                             2 1  1
\end{verbatim}
}
\smallskip
We construct a permutation with the corresponding wexbottoms and
wextops, and the corresponding nonzero values for $\Ca$ and $\Cb$,
that is, with $\Ca(5)=2$ and $\Cb(6)=\Cb(8)=1$.  First, the wexbottoms
are obtained by adding 1 to each descent bottom, and adjoining 1,
which is always a wexbottom.  The wextops are obtained by subtracting
1 from the descent tops, and adjoining $n$, which is always a wextop.
Thus, we get
\smallskip
\begin{verbatim}
        Wexbottoms:  (1) 2 4 5 7       Non-wexbots: 3 6 8 9

        Wextops:      1 5 6 8 (9)      Non-wextops: 2 3 4 7
\end{verbatim}
\smallskip
That is, $f$ is the word $12457$, $g$ is the word $3689$, 
$f'$ will be some permutation (to be determined) of the 
letters $15689$, and $g'$ will be some permutation (to be determined)
of the letters $2347$.

We construct the permutation in two parts, one for the weak
excedances, the other for the non-weak excedances.

Now, the definitions of $\Ca$ and $\Cb$ are such that $\Ca$ only
applies to pairs of weak excedances, and $\Cb$ only to pairs of
non-weak excedances.  We first construct the weak excedance part of
the permutation, by deciding where to place each of the wextops, in the
places given by the wexbots:

\begin{verbatim}
                              _ _ _ _ _
                              1 2 4 5 7

\end{verbatim}
We start from the right, in place 7, which has a 0 associated to it
(since $\rembr(7)=0$ in $\pi$).  We need to put there the smallest
number among the wextops that is at least as large as 7 (otherwise,
$\Ca(7)$ would exceed 0 in the resulting permutation).  This is the
number 8:

\begin{verbatim}
                                      8
\end{verbatim}
\bskip-1,
\begin{verbatim}
                              _ _ _ _ _
\end{verbatim}
\bskip-.3,
\begin{verbatim}
                              1 2 4 5 7
\end{verbatim}

\bskip.5, 

This leaves the wextops 1, 5, 6, 9.  The next place, 5, has a 2
associated to it (since $\rembr(5)=2$ in $\pi$), so we have to put
there a wextop that it is bigger than exactly two of the remaining
wextops that are at least as big as 5.  This forces us to make this 9
(and the two remaining wextops between 5 and 9 in size are 5 and 6).

We continue in this way until we have placed all the wextops, in such
a way that the values of $\ca$ for the remaining places are 0, since 5
is the only letter among the wexbottoms here with a nonzero right
embracing number in $\pi$:

\bskip.5,
\begin{verbatim}
                              1 6 5 9 8
\end{verbatim}
\bskip-1,
\begin{verbatim}
                              _ _ _ _ _
\end{verbatim}
\bskip-.3,
\begin{verbatim}
                              1 2 4 5 7
\end{verbatim}

\bskip.5,
The non-wex part is done in a similar way, but starting from the left,
and we get:

\bskip.5,
\begin{verbatim}
                              2 3 4 7
\end{verbatim}
\bskip-1,
\begin{verbatim}
                              _ _ _ _
\end{verbatim}
\bskip-.3,
\begin{verbatim}
                              3 6 8 9
\end{verbatim}
\bskip.5,
Observe that $\cb(6)=\cb(8)=1$ and $\cb(3)=\cb(9)=0$, as required.
Concatenating these two biwords, and sorting the columns to get the
bottom row in increasing order, the permutation we obtain is
$\Psi(215896374)=162593847$.

We now prove that the above procedure can always be carried out in the
way described.

\begin{theorem}\label{phi-thm}
  Let $\db'(\pi)$ be the set obtained from $\desbotset(\pi)$ by adding
  1 to each of its elements, and adjoining the letter 1.
  
  Let $\dt'(\pi)$ be the set obtained from $\destopset(\pi)$ by
  subtracting 1 from each of its elements, and adjoining the letter
  $n$.
  
  For a permutation $\tau$ let $\wb(\tau)$ be the set of weak
  excedance bottoms of $\tau$ and let $\wt(\tau)$ be the set of weak
  excedance tops of $\tau$.

  The map $\Psi$ described above is well defined, and has the
  following properties, where $\tau = \Psi(\pi)$:
\begin{itemize}
\item[(i)]%
  $\wb(\tau)=\db'(\pi)$,
\item[(ii)]%
  $\wt(\tau)=\dt'(\pi)$,
\item[(iii)]%
  $\ca(\tau)+\cb(\tau)=\rembr(\pi)$.
\end{itemize}
Moreover, $\Psi$ is a bijection.
\end{theorem}

\begin{proof}
  Recall that $\ca(i)=0$ unless $i$ is an excedance bottom, and that
  $\cb(i)=0$ unless $i$ is a non-excedance bottom.
  
  Let the letters of $\wb(\tau)$ be $b_1,b_2,\ldots,b_\ell$, ordered
  so that $ b_\ell<\cdots<b_2<b_1$. Look at the largest letter in
  $\wb(\tau)$, that is, $b_1$.  Suppose the embracing number of $b_1$
  in $\pi$ is $e_1$.  Then there are at least $e_1$ descent tops in
  $\pi$ that are larger than~$b_1$.  Thus, by the construction of
  $\wt(\tau)$ from the descent top set of $\pi$, there are at least
  $e_1+1$ elements $x$ in $\wt(\tau)$ such that $b_1\le x$.  So, we
  can find an element $t_1$ in $\wt(\tau)$ such that $\wt(\tau)$
  contains precisely $e_1$ elements $x$ satisfying $b_1\le x\le t_1$.
  Setting $\tau(b_1)=t_1$ guarantees that $\ca(b_1)=e_1$ in $\tau$.

  Look next at $b_2$, the second largest element in $\wb(\tau)$.
  Suppose its embracing number in $\pi$ is $e_2$.  There are then at
  least $e_2+1$ elements $x$ in $\wt(\tau)$ such that $b_2\le x$.
  However, one of these elements is $t_1$, which has already been
  placed to the right of place $b_2$ in $\tau$, and so $t_1$ cannot
  contribute to $\ca(b)$ in $\tau$.  But, $b_1+1$ is a descent bottom
  in $\pi$ and so its corresponding descent top, $d$, must be larger
  than $b_1+1$ and hence larger than $b_2$.  Thus, $b_2$ cannot be
  embraced by the descent $\ldots d(b_1+1)\ldots$ in $\pi$.  Hence,
  the embracing number of $b_2$ in $\pi$ can be at most one less than
  the number of elements $x$ in $\wt(\tau)$ satisfying $b_2\le x$.  We
  can therefore find an element $t_2\ne t_1$ in $\wt(\tau)$ such that
  precisely $e_2$ of the elements $x$ in $\wt(\tau)$ apart from $t_1$
  satisfy $b_2\le x\le t_2$.
  
  An analogous argument shows that the embracing number of $b_i$ in
  $\pi$ can be at most $N+1-i$, where $N$ is the number of elements
  $x$ in $\wt(\tau)$ with $b_i\le x$.  We can thus place each of the
  elements $t_i$ of $\wt(\tau)$ in $\tau$ so that $\ca(b_i)$ in $\tau$
  equals $\rembr(b_i)$ in $\pi$.
  
  In particular, each placement according to the above algorithm will
  result in the creation of a weak excedance.  Namely, clearly the
  $k$-th largest wexbottom is smaller than or equal to the $k$-th
  largest wextop.  Thus, by induction, since we consider the wexbots
  in decreasing order, the largest wextop unused at each stage of the
  algorithm is greater than or equal to the wexbottom being
  considered.
  
  To construct the subword of $\tau$ consisting of non-wextops, we
  proceed in a similar way, except that we start from the smallest
  non-wexbottom.  At each stage, for the non-wexbottom $b_i$ we find a
  non-wextop $d$ that satisfies $d<x<b_i$ for precisely $e$ elements
  $x$ among the remaining non-wextops, where $e=\rembr(b_i)$ in $\pi$.
  The argument showing that this is always possible, and that each
  placement results in a non-weak excedance, is analogous to the case
  of the weak excedance subword, and is omitted.
  
  To prove that $\Psi$ is a bijection, it suffices to show that it is
  injective, since it is a map from $\cls_n$ to itself.  Let
  $\sigma_1$ and $\sigma_2$ be two permutations with
  $\Psi(\sigma_1)=\Psi(\sigma_2)$.  From the definition of $\Psi$ it
  is clear that $\sigma_1$ and $\sigma_2$ must have the same descent
  tops and descent bottoms and also the same right embracing numbers
  for each letter.  It follows from the proof of Theorem~4 in \cite[p.
  249]{csz} that a permutation is uniquely determined by its sets of
  descent bottoms and tops, respectively, together with the right
  embracing numbers of its letters.  Thus, we must have
  $\sigma_1=\sigma_2$.
  
  In fact, the proof of Theorem~4 in \cite{csz} can be applied
  directly to our situation with trivial modifications, and yields a
  description of the inverse of $\Psi$.  In short, given the sets of
%
%
  weak excedance tops and weak excedance bottoms, respectively, of a
  permutation $\pi$, and the numbers $(\Ca(i)+\Cb(i))$, there is a
  unique permutation $\Psi^{-1}(\pi)$ with the corresponding descent
  tops and bottoms, respectively, and whose vector of right embracing
  numbers $\rembr(i)$ equals the vector of numbers $\Ca(i)+\Cb(i)$.
  The permutation $\Psi^{-1}(\pi)$ can be constructed in the exact
  same way as is done in the proof of Theorem~4 in \cite{csz}, except
  that the data we start with come from weak excedance tops and
  bottoms (and $\Ca$ and $\Cb$), instead of excedance tops and bottoms
  and the inversion numbers defined in \cite{csz}.
\end{proof}

Recall that $\wexbotsum(\pi)$ is the sum of all the wexbottoms in
$\pi$.  The following corollary of Theorem~\ref{phi-thm}  requires only
straightforward calculations.

\begin{corollary}\label{coro-wd}
\begin{eqnarray*}
\wextopsum\Psi(\pi) &=& \destopsum\pi + n - \des\pi,\\
\wexbotsum\Psi(\pi) &=& \desbotsum\pi + \des\pi + 1.
\end{eqnarray*}
\end{corollary}

We will use the following lemma
in our proofs of the equidistribution results
between our tableaux statistics and permutation statistics.

\begin{lemma}\label{lemma-a-w}
\begin{eqnarray}
\ac(\pi) &=& \ch n,2, - \ch n-\wex,2, + \wex - \wextopsum,\label{ac-wts}\\
\ad(\pi) &=& \wexbotsum - \ch\wex,2,.\label{ad-wbs}
\end{eqnarray}
\end{lemma}

\begin{proof}
  Equation~\pref{ac-wts} in the statement of the lemma is equivalent
  to
\begin{equation*}
\ac(\pi) + \wextopsum - \wex =  \ch n,2, - \ch n-\wex,2,.
\end{equation*}
We will show that the sum in the left-hand-side above
counts all pairs $(i,j)$,
with $1\le i<j\le n$, except those for which neither of $i$ and $j$ is
a weak excedance.

Recall that $\ac$ counts the pairs $(i,j)$ such that $j\le a_j<a_i<i$.
Each such pair can be described as consisting of a wextop $w$ in the
permutation, and a non-wextop that is larger than $w$ and to the right
of $w$.

We can interpret $(\wextopsum-\wex)$ as the sum, over all wextops, of
the size of the wextop, minus~1.  Counting this for each wextop $w$
can be done by counting all the letters in the permutation that are
strictly smaller than $w$.  Doing this for all wextops is equivalent
to counting all pairs of letters in the permutation that either
consist of two wextops, or a wextop and a non-wextop, where the wextop
is the larger of the two.

Therefore $\ac$ and $(\wextopsum-\wex)$ together count all pairs of
letters in the permutation, \emm except, those consisting of two
non-wextops.  (Observe that it is impossible to have a non-wextop $z$
and a wextop $w$ such that $z$ is left of $w$ and $z>w$.)  The total
number of pairs of letters in a permutation in $\cls_n$ is of course $\ch
n,2,$, and the number of pairs of non-wexbots is $\ch n-\wex,2,$,
which completes the proof.

Equation \pref{ad-wbs} in the statement of the lemma can be proved in
a similar manner.
\end{proof}

We can now prove the main results about the equidistribution implied
by the bijection $\Psi$.

\begin{thm}\label{equidist}
  Let $\sigma=\Psi(\pi)$.  We have
\begin{eqnarray}
 \des\pi &=& \wex\sigma-1,\label{dct-1}\\
 \caxb\pi &=& \aa(\sigma) + \ab(\sigma),\label{dct-2}\\
 \baxc\pi + \cxba\pi - \ch \des\pi,2,\label{dct-3} &=& \ac(\sigma),\\
 \bxca\pi &=& \ca(\sigma) + \cb(\sigma),\label{dct-4}\\
 \axcb\pi + \cbxa\pi - \ch \des\pi,2, &=& \ad(\sigma).\label{dct-5}
\end{eqnarray}
\end{thm}

\begin{proof}
  Equations~\pref{dct-1} and \pref{dct-4} in the statement of the
  theorem follow directly from Theorem~\ref{phi-thm}, since $\bxca\pi$
  is the sum of the right embracing numbers for all the letters in
  $\pi$.  We will prove \pref{dct-3} here; the proof of \pref{dct-5}
  is analogous and is omitted.  Having done this,
  Equation~\pref{dct-2} follows from the other four identites in the
  present theorem, together with Lemmas~\ref{lemma-corteel}
  and~\ref{patt-const} and routine calculations.
  
  To prove Equation~\pref{dct-3}, observe that
\begin{equation*}
\baxc\pi + \cxba\pi  = n\cdot\des\pi - \destopsum\pi.
\end{equation*}
This is because $\baxc\pi+\cxba\pi$ counts
the letters in $\pi$ larger than the descent top $b$
for each descent $\ldots ba\ldots$ in $\pi$.
According to Corollary~\ref{coro-wd}, the right-hand-side
in the equation above
can be rewritten as follows:
\begin{equation*}
  n\cdot\des\pi - \destopsum = n\cdot\des\pi - \wextopsum\sigma + n -
  \des\pi.
\end{equation*}
By Lemma~\ref{lemma-a-w}, this is equal to
\begin{equation*}
n\cdot\des\pi+\left(\ac(\sigma)-\wex\sigma-\ch n,2,+\ch n-\wex\sigma,2,\right)+n-\des\pi,
\end{equation*}
which, in turn, is equal to
\begin{equation*}
\ac(\sigma)+ n\cdot\des\pi-(\des\pi+1)-\ch n,2,+
\ch n-(\des\pi+1),2,+n-\des\pi.
\end{equation*}
To show that this last expression is equal to
 $\ac(\sigma)+\ch\des\pi,2,$ is straightforward.
\end{proof}

Note that Equations \pref{dct-2} and \pref{dct-3} in Theorem
\ref{equidist} together imply that
\begin{equation*}
  \caxb\pi + \baxc\pi + \cxba\pi  - \ch \des,2,= 
  \aa(\sigma) + \ab(\sigma) + \ac(\sigma).
\end{equation*}
This last equation,
 together with Theorem~\ref{equidistribution1}, leads to the
following corollary.

\begin{corollary}\label{equidistribution2}
 Let $T(k,a,b,c)$ be the set of permutation tableaux with $k$
  rows and $(n-k)$ columns, which are filled with precisely
 $a$ 0's, $b$ 1's and $c$ 2's.  Let
  $P(k,a,b,c)$  be the set of all permutations $\pi \in \cls_n$, such that
\begin{itemize}
\item%
 $k-1 = \des(\pi)$,
\item%
  $a = \left[\caxb + \baxc +\cxba\right]\pi -\ch\des\pi,2,$,
\item%
  $b = \bxca\pi + n-1-\des\pi$,
\item%
  $c =\left[\axcb + \cbxa\right]\pi - \ch \des\pi,2,$.
\end{itemize}
Then $|T(d,a,b,c)| = |P(d,a,b,c)|$.
\end{corollary}

\section{Enumeration results}\label{Enumeration}

One nice application of permutation tableaux is that they facilitate
enumeration of permutations according to various statistics.  This is
because permutation tableaux satisfy a rather simple recurrence, which
we now explain.

\begin{figure}[h]
  \centerline{\epsfig{figure=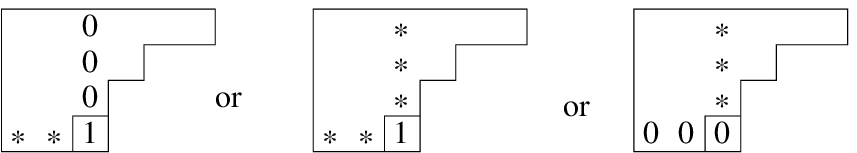}}
\caption{Recurrence for $F_{\lambda}(p, q)$}
\label{Recurrence}
\end{figure}

Fix a partition $\lambda = (\lambda_1, \dots , \lambda_k)$.  Let
$F_{\lambda} (p, q)$ be the polynomial in $p$ and $q$ such that the
coefficient of $p^s q^t$ is the number of valid fillings of the Young
diagram $Y_{\lambda}$ which contain $s$ $0$'s and $t$ $1$'s.  As
Figure \ref{Recurrence} illustrates, there is a simple recurrence for
$F_{\lambda}(p, q)$.  (Note, however, that this is {\it not} the 
same as the recurrence given for $\Le$-diagrams in \cite{Williams}.)

Explicitly, any valid filling of $\lambda$ is obtained in one of the
following ways:
\begin{itemize}
\item inserting a column whose bottom entry is $1$ and whose other
  entries are $0$ after the $(\lambda_k - 1)$st column of a valid
  filling of $(\lambda_1 - 1, \lambda_2 - 1, \dots , \lambda_k - 1)$;
\item adding a $1$ to the end of the bottom row of a valid filling of
  the shape $(\lambda_1, \lambda_2, \dots , \lambda_{k-1}, \lambda_{k}
  - 1)$;
\item adding an all-zero row of length $\lambda_{k-1}$ to a valid
  filling of $(\lambda_1, \dots , \lambda_{k-1})$.
\end{itemize}

Thus, we have the following recurrence.

\begin{prop}\label{recur}
\begin{multline*}
F_{\lambda}(p,q) = p^{k-1}q F_{(\lambda_1 - 1, \lambda_2 - 1, \dots,
\lambda_k -1)}(p,q) + q F_{(\lambda_1, \lambda_2,\dots, \lambda_{k-1},
\lambda_{k} -1)}(p,q)\\ + p^{\lambda_k} F_{(\lambda_1, \dots,
\lambda_{k-1})}(p,q).
\end{multline*}
\end{prop}

It is straightforward to compute $F_{\lambda}(p,q)$ when
$k$ (the number of rows of $\lambda$) 
is small.  Here are the first two formulas.

\begin{prop}
\begin{align*}
F_{(\lambda_1)}(p,q) &= q^{\lambda_1}. \\
F_{(\lambda_1, \lambda_2)}(p,q) &= -q^{\lambda_1 - 1}p^{\lambda_2 + 1} +
q^{\lambda_1 - 1}[2]_{p,q}^{\lambda_2 + 1}.
\end{align*}
\end{prop}

In the above expression, $[2]_{p,q}$ is the $p,q$-analog of $2$.  
Recall that the $p,q$-analog of the number $n$ is 
$p^{n-1} + p^{n-2}q + p^{n-3}q^2 + \dots + q^{n-1}$, 
denoted $[n]_{p,q}$.

Let $D_{k,n}(p,q,r):= \sum_{\lambda}
F_{\lambda}(p,q)r^{k(n-k)-|\lambda|}$, where $\lambda$ ranges over all
partitions contained in a $k \times (n-k)$ rectangle.  By Theorem
\ref{equidistribution1}, $D_{k,n}(p,q,r)$ enumerates permutations
according to the number of weak excedances, several kinds of
alignments, and crossings.  And by Corollary \ref{equidistribution2},
$D_{k,n}(p,q,r)$ enumerates permutations according to the number of
descents and occurrences of various generalized patterns.  Therefore
it would be nice to get an explicit expression for $D_{k,n}(p,q,r)$,
for example by solving explicitly for $F_{\lambda}(p,q)$ and then by
summing over partitions $\lambda$ contained in a $k \times (n-k)$
rectangle.

For fixed small $k$, it is not too difficult to compute the generating
function $D_{k}(p,q,r,x):= \sum_n D_{k,n}(p,q,r) x^n$.  Here are the
first few formulas.  Note that it is easy to determine what the
denominator should be for $D_k(p,q,r,x)$, but the numerator is
significantly more complicated.

\begin{prop}
\begin{align*}
&D_1(p,q,r,x) = \frac{x}{1-qx}\\
&D_2(p,q,r,x) = \frac{x^2}{(1-pqx)(1-qrx)(1-q[2]_{p,q}x)}\\
&D_3(p,q,r,x) =\\ 
&\phantom{aaa}\frac{x^3(1+pq^2x - p^3q^2rx^2-2p^2q^3rx^2 - pq^4rx^2)}
{(1-p^2qx)(1-pqrx)(1-qr^2x)(1-pq[2]_{p,q}x)(1-qr[2]_{p,q}x)(1-q[3]_{p,q}x)}
\end{align*}
\end{prop}
\bigskip

One can derive these formulas by either using the methodology outlined
above (i.e.\ by summing $F_{\lambda}(p,q)$), or else by translating
the problem of enumerating permutation tableaux into a problem about
enumerating certain weighted lattice paths, and then by 
enumerating these lattice paths.  In order to sketch
the latter method, let us define a ${\it bad}$ zero in a permutation 
tableau to be a $0$ which lies directly underneath some $1$.  Note that
if some column $C$ in a permutation tableau $\T$
contains a bad zero in the $r$th row, then every column 
to the left of $C$ must also contain a zero in the $r$th row.

In the lattice path method for enumeration of permutation tableaux,
we associate to each permutation tableau $\T$ a weighted
lattice path $L = \{L_i\}_{i=1}^n$ 
consisting of $n$ steps in the plane, which must be of the following types:
$(1,1)$ (a northeast step), 
$(1,0)$ (an east step), 
and $(1, -j)$, where $1 \leq j \leq k-1$ (a southeast step).  Each 
step $L_i$ in the lattice path  represents the step
$P_i$ 
in the partition path $\{P_i\}_{i=1}^n$.  (Recall that the partition path
follows the shape of the partition $Y_{\lambda}$ and travels
from the northeast corner to the 
southwest corner of the $k \times (n-k)$ rectangle containing
$\T$.)  
The steps $(1,1)$ in $L$ correspond to 
vertical steps in the partition path, and have weight $x$. 
A step $(1,0)$ in $L$ corresponds to a horizontal 
step in the partition path such that the corresponding column $C$
of $\T$ does not introduce any bad 
zeros (except those that were forced by bad zeros to the right of $C$).
Such a step has weight 
$p^a q^b r^c x$, where $a$, $b$, and $c$ are the numbers of $0$'s, 
$1$'s, and $2$'s, respectively, in column $C$. (Note that $a+b+c=k$.)
Finally, a step $(1, -j)$ in $L$ corresponds to a horizontal 
step in the partition path such that the corresponding column $C$ of 
$\T$ introduces exactly $j$ new 
bad zeros (that were not forced by bad zeros in columns 
to the right of $C$).
As before, such a step has weight 
$p^a q^b r^c x$, where $a$, $b$, and $c$ are the numbers of $0$'s, 
$1$'s, and $2$'s, respectively, in column $C$.
Observe that the height of any point in the lattice path $L$
is equal to the number of boxes of the corresponding column of $\T$
which can be filled with either a $0$ or a $1$.
By associating weighted lattice paths to permutation tableaux in this
way, we can facilitate computation of the generating functions 
$D_k(p,q,r,x)$ for small $k$.

Now we will give complete results about a certain specialization
of $D_{k,n}(p,q,r)$.  
Let $E_{k,n}(q):= D_{k,n}(1,q,1)$.  An explicit formula for 
$E_{k,n}(q)$ was found in \cite{Williams}; the proof 
utilized  
a recurrence similar to that in Proposition \ref{recur}.

\begin{theorem}[Williams \cite{Williams}]
\begin{equation*}
{E}_{k,n}(q) =
          q^{n-k^2} \sum_{i=0}^{k-1} (-1)^i [k-i]^n q^{ki-k}
       \left( {n \choose i} q^{k-i} + {n \choose i-1}\right).
\end{equation*} \label{formula}
\end{theorem}

In the above formula, the notation $[k-i]$ refers to the 
$q$-analog of the number $k-i$, that is, 
$1+q+q^2+\dots + q^{k-i-1}$.

The polynomials above have many nice properties.
It was observed in \cite{Williams} that if one 
renormalizes $E_{k,n}(q)$ by defining
$\hat{E}_{k,n}(q):= q^{k-n}E_{k,n}(q)$, then 
$\hat{E}_{k,n}(q)$ is a new $q$-analog of the Eulerian numbers
(distinct from Carlitz' classical
$q$-analog of the Eulerian numbers \cite{Carlitz}).
Furthermore, $\hat{E}_{k,n}(q)$
specializes at $q=-1, 0, 1$ to 
the binomial coefficients, the Narayana numbers, and the Eulerian numbers.
Additionally, $\hat{E}_{k,n}(q) = \hat{E}_{n+1-k,n}(q)$.
It was shown more recently by
Corteel \cite{Corteel} that the polynomials 
$\hat{E}_{k,n}(q)$ naturally relate to the ASEP model in statistical physics.

Theorem \ref{formula} together with Corollary \ref{equidistribution2}
implies the following result.
\begin{corollary}
The number of permutations in $\cls_n$ with $k-1$ descents and $m$
occurrences of the pattern $(2-31)$ is equal to the coefficient of
$q^m$ in
\begin{equation*}
\hat{E}_{k,n}(q)=
          q^{-k^2} \sum_{i=0}^{k-1} (-1)^i [k-i]^n q^{ki}
       \left( {n \choose i} q^{k-i} + {n \choose i-1}\right).
\end{equation*}
\end{corollary}
This result was first conjectured by the authors of this paper, and
first proved by Corteel \cite{Corteel}.  The formula
$\hat{E}_{k,n}(q)$ is the first known polynomial expression which
gives the {\em complete} distribution of a permutation pattern of
length greater than 2 (the two cases of length 2 correspond to the
Eulerian numbers and the coefficients of $[n]!$, respectively).

The generating function for the polynomials $\hat{E}_{k,n}(q)$
has been expressed in two ways: as a formal power series
and as a continued fraction.
That is, it can be shown
\cite{Williams} that $\hat{E}(q,x,y):=\sum_{n,k}\hat{E}_{k,n}(q)y^k
x^n$ is equal to
\begin{equation*}
\sum_{i=0}^{\infty} \frac{y^i (q^{2i+1} - y)}
{q^{i^2+i+1} (q^i - q^{i+1}[i]x + [i]xy)}.
\end{equation*}

Additionally, Corteel \cite{Corteel} used 
results of Clark, Steingr\'{\i}msson, and Zeng
\cite{csz} to show the following:
\begin{theorem}[Corteel \cite{Corteel}]
\begin{equation*}
\hat{E}(q,x,y)=
\cfrac{1}{1-b_0x-\cfrac{\lambda_1x^2}{1-b_1x-\cfrac{\lambda_2x^2}
{1-b_2x-\cfrac{\lambda_3x^2}{\ddots}}}},
\end{equation*}
\null
\bigskip
where $b_n=y[n+1]_q+[n]_q$, $\lambda_n=y[n]_q^2$,
and $[n]_q=1+q+\cdots +q^{n-1}$.
\label{big}
\end{theorem}

\subsection{The Euler-Mahonian distribution and 
Carlitz' $q$-analog of the Eulerian numbers}

Recall that the {\it Euler-Mahonian} distribution is the joint
distribution of the number of descents and the {\it major index}
for permutations in $S_n$.
The {\it major index} is the sum of the places of the descents in 
a permutation.

The generating function for this joint distribution is 
given by Carlitz' $q$-Eulerian polynomials $B_{n,k}(q)$ \cite{Carlitz}, 
which one can define by
\begin{equation}
B_{n,k}(q) = \sum_{\pi} q^{\expmaj(\pi)-{k \choose 2}},
\end{equation}
where the sum is over permutations in $S_n$ which have $k-1$ descents.
Note that one subtracts ${k \choose 2}$ from the exponent because 
when the number of descents of $\pi$ is $k-1$, 
the quantity $\expmaj(\pi)$ is at least ${k \choose 2}$.

Analogous to the Eulerian numbers, the coefficients of the $q$-Eulerian
polynomial satisfy the recurrence \cite{Carlitz}
\begin{equation}\label{Carlitz-recurrence}
B_{n,k}(q) = [k+1]B_{n-1,k}(q)+q^k [n-k] B_{n-1,k-1}(q),
\end{equation}
subject to the initial conditions $B_{0,k}(q)=1$ for $k=0$, and
$B_{0,k}(q)=0$ otherwise.

We will now show that Carlitz' $q$-analog is simply 
a specialization of the 
polynomial $D_{k,n}(p,q,r)$.

\begin{prop}\label{specialization}
With the notation above, we have that 
$D_{k,n}(p,1,1) = B_{n,k}(p)$.
\end{prop}

\begin{proof}
We will prove this by 
showing that if
$(\mathrm{rows},\mathrm{zeros})$ is the bistatistic counting rows
and 0's in tableaux, then $(\mathrm{rows},\mathrm{zeros})$ has the
same distribution as the pair $(\des+1, \maj-\ch\des+1,2,)$.

Note first that the statistic $a$ in
Corollary~\ref{equidistribution2}, when stripped of $\ch\des,2,$, has
the same distribution as the statistic
\begin{equation}\label{mak}
\axcb+\cbxa+\bxca.
\end{equation}
This is because the statistic in \pref{mak} is obtained by taking the
\emm reverse complement, of the statistic $\caxb+\baxc+\cxba$, that
is, by reversing each of the patterns and then replacing each letter
$i$ by $4-i$.  Doing the same with each permutation in $\cls_n$ (with
4 replaced by $n+1$) is a bijection from $\cls_n$ to itself, and this
bijection clearly proves the equidistribution of $\caxb+\baxc+\cxba$
with $\axcb+\cbxa+\bxca$, even when each statistic is taken jointly
with the number of descents (which is invariant under reverse
complement).  The statistic $\axcb+\cbxa+\bxca+\des$ is equal to the
statistic $\mak$, as pointed out in \cite{bast}, and it was shown by
Foata and Zeilberger \cite{foata-zeil} that $(\des,\mak)$ has the
Euler-Mahonian distribution, that is, the same distribution as
$(\des,\maj)$.\footnote{Actually, the statistic $\axcb+\cbxa+\bxca$ is
a slight variation on $\mak$, as defined by Foata and Zeilberger, but
is easily seen to have the same distribution when taken together with
the number of descents. Foata and Zeilberger's $\mak$ is actually
equal to the statistic called $\makl$ in \cite[Table 1]{bast}} It
follows that $(\des,\maj-\ch\des+1,2,)=(\des,\maj-\des-\ch\des,2,)$
has the same distribution as $(\des,\mak-\des-\ch\des,2,)$, which, in
turn, has the same distribution as
$(\des,\axcb+\cbxa+\bxca-\ch\des,2,)$ .  Hence, by Corollary
\ref{equidistribution2}, $(\des+1,\maj-\ch\des+1,2,)$ has the same
distribution as $(\mathrm{rows},\mathrm{zeros})$.
  
Therefore $D_{k,n}(p, 1, 1)
=B_{n,k}(p)$.
\end{proof}

Thus the polynomials $D_{k,n}(p,q,r)$ generalize both the classical
$q$-analog of the Eulerian numbers and the new $q$-analog of the
Eulerian numbers found in \cite{Williams}.

\section{Open problems and other remarks}

When this section of the paper was first written, it contained a 
collection of open problems that we thought were worth studying.
We are happy to report that in the nine months following
our posting of this paper on the 
electronic arXiv, 
there was a great deal of progress on our open problems
by Alexander Burstein \cite{Alex}, Sylvie Corteel \cite{Corteel0},
Niklas Eriksen \cite{Niklas},
Astrid Reifegerste \cite{Astrid}, and Xavier Viennot \cite{Xavier}.
We will list the full set of open problems below, with remarks at the 
end about the 
progress that has been made on them.

\begin{enumerate}

\item \label{1}
Find an explicit expression for $D_{k,n}(p,q,r)$.

\item  \label{2}
Can one prove Proposition \ref{specialization} by 
checking the analogous recurrence (\ref{Carlitz-recurrence})
for permutation tableaux?

\item \label{3}
We say that a 1 in a permutation tableau is \emm essential, if it is
  the topmost one in its column or the leftmost 1 in its row.  A
  tableau is determined by its essential 1's: all the other 1's
  are determined by these, because of condition (2)
  in the definition of a permutation tableau.  What do the essential
  1's correspond to in the corresponding permutation?  

 We conjecture---based on experimental evidence for $n$ up to $9$
 ---that the distribution of permutation tableaux according to the
 number of essential 1's is equal to that of permutations according to
 $(n-c)$, where $n$ is the length of the permutations and $c$ the
 number of cycles when each permutation is written in standard cycle
 form.  This distribution is the same as that for Left-to-Right
 minima.  Moreover, we conjecture that the joint distribution of
 tableaux according to the number of rows and the number of essential
 1's equals that of permutations according to $(n-1-\des)$ and
 $(n-\mathrm{LR})$, where des is the number of descents, and LR the
 number of Left-to-Right-minima.  The bistatistic
 $(\des+1,\mathrm{LR})$, in turn, has the same distribution as the
 number of weak excedances and the number of cycles of a permutation,
 when written in standard cycle form.

\item \label{4}
 The number of 0's in a tableau corresponds to the total number of
  occurrences of the patterns $\cxba$, $\baxc$ and $\caxb$.  It is
  easy to see that these patterns have the same distributions as
  $\axcb$, $\cbxa$ and $\bxca$, respectively.  To prove this, simply
  reverse each permutation in $\cls_n$.  Can we partition the 0's in a
  tableau into two sets, one corresponding to occurrences of
  $\cxba+\baxc$ and the other to occurrences of $\caxb$?  Observe that
  the first one of these sets would correspond to descent tops and the
  second one to left embracings.  Thus, these sets would be symmetric
  counterparts of 2's and 1's respectively, although this symmetry is
  not transparent in the tableaux.

\item \label{5}
  The reflection of a permutation tableau $\T$ in its
  north-west/south-east diagonal yields a permutation tableau if and
  only if $\T$ has a 1 in each row.  That is equivalent to the
  associated permutation being fixed point free.  Which permutation is
  associated to the reflected tableau (that tableau is also fixed
  point free because it has a 1 in each row)?

\item \label{6}
A permutation tableau $T$ must have at least one 1 in each column.  If
it has only this minimum number of 1's, then the corresponding
permutation, that is, $\Psi^{-1}(\Phi(T))$, has no occurrences of the
pattern $\bxca$.  It has been shown (see \cite{claesson}), that
permutations avoiding this pattern are enumerated by the Catalan
numbers.  Is there a bijection from these tableaux to any of the well
known objects enumerated by Catalan numbers, such as Dyck paths?

\item \label{7}
Find a better description of the bijection $\Psi$, and of the 
composition of maps $\Psi^{-1} \circ \Phi$.
\end{enumerate}

As mentioned earlier, much progress on these problems has been made by
Alexander Burstein \cite{Alex}, Sylvie Corteel \cite{Corteel0}, 
Niklas Eriksen \cite{Niklas},
Astrid Reifegerste \cite{Astrid}, and Xavier Viennot \cite{Xavier}.

Independently, Alexander Burstein \cite{Alex} and Niklas Eriksen
\cite{Niklas} have solved open problems \ref{3}, \ref{5}, \ref{6}.
Additionally, 
Sylvie Corteel \cite{Corteel0} and Astrid Reifegerste \cite{Astrid} have 
(independently) solved problem \ref{6}.
Xavier Viennot \cite{Xavier} has found a new bijection from permutation tableaux
to permutations which answers the open problems \ref{2} and also \ref{6}.
Additionally, it seems to be related to the recent work
of the second author \cite{Williams3}
on the asymmetric exclusion process.

\section{Acknowledgements}

We would like to thank the two referees for a
careful reading of the paper and many comments, which resulted in a
significant improvement of the presentation.  Additionally, 
the second author would like to think Alex Postnikov for interesting
discussions.

\end{document}